\documentclass[11pt]{article}

\usepackage{amsmath}
\usepackage{amssymb}

\newtheorem{theorem}{Theorem}
\newtheorem{definition}{Definition}
\newtheorem{lemma}{Lemma}
\newtheorem{proposition}{Proposition}

\newcommand{\all}{\forall}
\newcommand{\ext}{\exists}

\newcommand{\sbi}[1]{\Sigma^b_{#1}}

\newcommand{\sBi}[1]{\Sigma^B_{#1}}
\newcommand{\pBi}[1]{\Pi^B_{#1}}
\newcommand{\dBi}[1]{\Delta^B_{#1}}

\newcommand{\ltwo}{\mathcal{L}_A^2}
\newcommand{\pairs}[1]{\langle #1\rangle}
\newcommand{\ptran}[1]{\Vert #1\Vert}
\newcommand{\bv}[1]{\ensuremath{[\![#1]\!]} }

\newcommand{\circuits}{\mathbb{C}}
\newcommand{\boolec}{\mathbb{B}}
\newcommand{\boola}{\mathbb{B}_A}
\newcommand{\boolef}{\mathbb{B}_{EF}}
\newcommand{\booler}{\mathbb{B}_R}
\newcommand{\lpv}{\mathcal{L}^2_{PV}}

\newcommand{\apc}{APC_1}

\newcommand{\pv}{\mathbf{PV}}
\newcommand{\vp}{\mathbf{VP}}
\newcommand{\vi}[1]{\mathbf{V}^{#1}}
\newcommand{\rcircuit}{\mathbb{RC}(\bar{p})}

\title{Developing Takeuti-Yasumoto forcing}

\author{Satoru Kuroda\\Gunma Prefectural Women's University}

\begin{document}

\maketitle

\section{Introduction}
\label{sec:introduction}

In late 90's two consecutive papers by G. Takeuti and Y. Yasumoto \cite{ty}, \cite{ty2}
were published. These papers aimed at applying Boolean valued model constructions
for bounded arithmetic and relating them to separation problems of complexity classes.

Prior to their work, forcing was used to construct models of bounded arithmetic
in different contexts. The first application of forcing in bounded arithmetic
was done by Paris and Wilkie, who proved that the theory $IE_1$ extended by a single
predicate symbol $R$ does not prove that $R$ is not a bijection from $n+1$ to $n$.
Their construction was later extended to $I\Delta_0$ by M. Ajtai and to
Buss' theory by S. Riis.

Besides these results, J. Kraj\'\i\v cek \cite{jan} gave different type of
forcing construction for models of bounded arithmetic $PV$ and $V^1$. Kraj\'\i\v cek's
motivation for these results is to provide nonstandard models of weak theories
which satisfies unproven separation of complexity classes.

Then Takeuti and Yasumoto tried to give a comprehensive theory for forcing
in nonstandard models of bounded arithmetic. Their main motivation was to
relate separations of complexity classes in the standard world to generic models
constructed from nonstandard models.

Although forcing type arguments seems very useful in bounded arithmetic
just like in other branches of mathematical logic, there are few connections
between the above mentioned results.

In this paper, we rearrange Takeuti-Yasumoto type forcing argument in a frame work
which is different from the one they adopted. Namely we rework on the construction
of generic models in two sort language in Cook and Nguyen \cite{cn}.
It seems that some of the proofs in \cite{ty} and \cite{ty2} have flaws and
we give a correct proofs of some of their results.

We will also show that Kraj\'\i\v cek's results can be obtained as Takeuti-Yasumoto
forcing construction. Furthermore, we will consider the problem of violating
or satisfying surjective weak pigeonhole principles for polynomial time functions
in generic models. Actually, we show that under an assumption on propositional
logic in the ground model, we can construct a generic model violating $dWPHP(PV)$.
On the other hand, we can construct generic models for $dWPHP(PV)$ by extending
the base Boolean algebra to allow random inputs.

Our results suggest that Takeuti-Yasumoto forcing works as a general framework
for forcing construtions in bounded arithmetic.

\section{Preliminaries}
\label{sec:preliminaries}

We will work on two-sort bounded arithmetic which was developed by Cook and Nguyen \cite{cn}.
The language of two-sort bounded arithmetic comprises two sorts of variables;
number variables $x,y,z,\ldots$ and string variables $X,Y,Z,\ldots$.
The language $\ltwo$ has a constant symbol $0$, function symbols $s(x)$, $x+y$,
$x\cdot y$, $|X|$ and a relation symbol $x\leq y$ where $|X|$ denotes the length
of string $X$. We use either expressions $X(i)$ and $i\in X$ to denote that
the $i$th bit of $X$ is $1$.

A bounded number quantifier is of the form $\all x<t$ or $\ext x<t$. A bounded string
quantifier is of the form $\all X<t$ or $\ext X<t$ whose intended meanings are
$$
\all X(|X|<t\rightarrow\cdots)\mbox{ and }\ext X(|X|<t\land\cdots)
$$
respectively. $\sBi{0}$ is the set of formulas which contains only bounded number
quantifiers. $\sBi{1}$ is the set of formulas which contains bounded number quantifiers,
positive occurrences of bounded existential string quantifiers and
negative occurrences of bounded existential string quantifiers.

In this paper we treat theories for PTIME  and their extensions. In particular
the following three theories for PTIME  are considered.

$\vp$ is the $\ltwo$ theory whose axioms are $BASIC_2$ and a single axiom $MCV$
given as
$$
MCV\equiv\all a\;\all C,E\;\ext Y\;\delta_{MCV}(a,C,E,Y)
$$
where
$$
\begin{array}{ll}
\delta_{MCV}(a,C,E,Y)\equiv\\
\neg Y(0)\land Y(1)\land \all x<a\; x\geq 2\rightarrow\\
Y(x)\leftrightarrow[(C(x)\land\all y<x\;(E(y,x)\rightarrow Y(y))\lor
                    (\neg C(x)\land\ext y<x\;(E(y,x)\land Y(y)))].
\end{array}
$$

The second theory $\pv$ is defined over the language $\lpv$ which extends $\ltwo$
by function symbols for Cobham's function algebra for PTIME. Then $\pv$ consists of
$BASIC_2$ and defining axioms for such functions.

A seemingly stronger theory $\vi{1}$ is an $\ltwo$-theory which comprises $BASIC_2$
together with
$$
\sBi{1}\mbox{-COMP : }\ext Y<a\;\all x<a\;(Y(x)\leftrightarrow\varphi(x)),\ \varphi(x)\in\sBi{1}.
$$

It is known that all these theories corresponds to the class PTIME in the sense that
\begin{theorem}[Cook-Nguyen \cite{cn}]
For $T=\vp$, $\pv$ or $\vi{1}$, a function is $\sBi{1}$-definable in $T$ if and only if
it is computable in PTIME.
\end{theorem}

We will work with the circuit models for PTIME. In the axiom $MCV$, monotone circuits are
coded by pairs of the form $(G,E)$ of strings such that
$$
Circuit_n(G,E)\leftrightarrow E\subseteq G\times G.
$$

The following fact is folklore and will be used elsewhere in this paper.
\begin{proposition}
Let $F(X)$ be a $PV$ function. Then $\pv$ proves the followng.
$$
\all n\;\ext\pairs{C_0,\ldots,C_{t-1}}\;\all i<t(C_i(X)=1\leftrightarrow F(X)(i))
$$
where $t$ is the bounding term for $F$.
\end{proposition}

Next we review basic notions of Boolean valued models in \cite{ty} in terms of
the two-sort bounded arithmetic.

Let $(M_0,M)$ be a structure of some two-sort language of bounded arithmetic.
Throughout the paper we will concentrate on nonstandard countable model
$(M_0,M)\models VP+\neg Exp$.

\noindent
{\bf Example.}
The original Takeuti-Yasumoto forcing starts with the ground model which is defined in the
following manner. Let $M\models Th(\mathbb{N})$ be a countable nonstandard model and
fix $n\in M\setminus\omega$. Define
$$
M^*=\{x\in M : x\leq \underbrace{n\#\cdots\#n}_{k\mbox{ times}} \mbox{ for some }k\in\omega\}
$$
and
$$
M_0=\{|x| : x\in M^*\}.
$$
We regard $(M_0,M^*)$ as a two-sort structure by identifying each element $x\in M^*$
with its binary representation.
It is easy to see that $(M_0,M^*)\models \vi{\infty}$.

Let $n_0=|n|$ and $\bar{p}=p_0,\ldots,p_{n_0-1}$ be the list of propositional variables
coded by elements in $M_0$.
We define $\mathbb{C}$ as a set of Boolean formulas over variables from $\bar{p}$
coded in $M^*$.
The precise definition of $\mathbb{C}$ can be found in \cite{ty}.

$\mathbb{C}$ can be regarded as a Boolean algebra with respect to either one of
the following two partial orders:
$$
\begin{array}{ll}
  C\leq_A C'& \Leftrightarrow\all X (|X|=n\rightarrow eval(C,X)\leq eval(C',X)),\\
  C\leq_{EF} C'&\Leftrightarrow\ext P\; Prf_{Ef}(C\rightarrow C',P).\\
\end{array}
$$
Note that these two partial orders are identical only if Extended Frege is super.
In the following definitions, the partial order $\leq$ on $\mathbb{C}$ represents
either $\leq_A$ or $\leq_{EF}$.

Define $\boola=\circuits/=_A$ and $\boolef=\circuits/=_{EF}$. We omit the subscript
and denote either Boolean algebra by $\boolec$ when there is no fear of confusion.

A set $I\subseteq\mathbb{B}$ is an ideal if $0\in I$, $1\not\in I$ it is closed
under $\lor$ and lower closed with respect to the partial order.
An ideal $I$ is $M_0$-complete if
$$
\all n\in M_0\;\all X:a\rightarrow\mathbb{B}\;\all i<a\;
X(i)\in I\Rightarrow\bigvee_{i<a}X(i)\in I.
$$

A set $F\subseteq\mathbb{B}$ is a filter if $0\not\in I$, $1\in I$ it is closed
under $\land$ and upper closed　with respect to the partial order.

A set $D\subseteq\mathbb{B}$ is dense over an $M_0$-complete ideal $I$
if for any $X\in\mathbb{B}\setminus I$ there is $X'\in \mathbb{B}\setminus I$
such that $X'\leq X$. $D$ is definable if there exists a formula $\varphi$ such that
$$
D=\{X\in\mathbb{B} : M\models\varphi(X)\}.
$$

A filter $G\subseteq\mathbb{B}$ is $\mathcal{M}$-generic
if $(D\setminus I)\cap G\neq\emptyset$ whenever $D$ is dense over $I$ and definable.
Remark that an $\mathcal{M}$-generic $G$ is not definable in $(M_0,M)$.

We define the generic model analogous to that in set theory. First define
$$
M^{\mathbb{B}}=\{X\in M : X:a\rightarrow\mathbb{B}\mbox{ for some }a\in M_0\}.
$$
For a $\mathcal{M}$-generic $G$ over an $M_0$-complete ideal $I$ and $X:a\rightarrow\mathbb{B}$ we define
$$
i_G(X)=\{y<a : X(y)\in G\}.
$$
Finally we define
$$
M[G]=\{i_G(X) : X\in\mathbb{B}\}.
$$
Since the length of any $i_G(X)$ is bounded by some element in $M_0$, we can regard the pair $(M_0,M[G])$
as a two-sort structure with a natural interpretation.

First we will show that $(M_0,M[G])$ is a model of the base theory $V^0$.

\begin{definition}
Let $\varphi(\bar{x},\bar{X})\in\sBi{0}$ and
$$
\ptran{\varphi(\bar{x},\bar{X})}_{\bar{m},\bar{n}}(\bar{p}_0,\ldots,\bar{p}_{n-1})
$$
be its propositional translation where $\bar{p}_i$ corresponds to the variable $X_i$.
For $\bar{a}\in M_0$ and $\bar{A}\in M^{\boolec}$ with $A_i:b_i\rightarrow\boolec$,
we define
$$
\bv{\varphi(\bar{a},\bar{A})}
\ptran{\varphi(\bar{x},\bar{X})}_{\bar{a},\bar{b}}(X_0,\ldots,X_{n-1}).
$$
\end{definition}

\begin{theorem}[Forcing Theorem for $\sBi{0}$ formulas]
Let $\mathbb{B}$ be a Boolean algebra in $M^*$ and suppose that it admits a $\sBi{0}$-translation.
Let $\varphi(\bar{x},\bar{X})$ be a $\sBi{0}$ formula with parameters as indicated.
Then for $\bar{a}\in M_0$ and $\bar{A}\in M^{\mathbb{B}}$
$$
(M_0,M[G])\models\varphi(\bar{a},\bar{i_G(X)})\Leftrightarrow\bv{\varphi}\in G
$$
for a $\mathcal{M}$-generic $G$ over an $M_0$-complete ideal $I$
\end{theorem}

\begin{theorem}
Let $\mathbb{B}$ be a Boolean algebra in $M^*$ and suppose that it admits a $\sBi{0}$-translation.
If $G\subseteq\mathbb{B}$ is a $\mathcal{M}$-generic over an $M_0$-complete ideal $I\subseteq\mathbb{B}$ then
$$
(M_0,M[G])\models V^0.
$$
\end{theorem}

For the Boolean algebra $\boolec$, we have a stronger forcing theorem, that is
$\sbi{0}$ formula may contain $PV$-functions.

\begin{theorem}[Cobham, Cook-Nguyen]
A string function is in $FP$ if and only if it is obtained by $AC^0$ functions
by finitely many applications of composition and limited recursion.
\end{theorem}

\begin{lemma}\label{lem:bound}
  If $F(\bar{x},\bar{X})\in\lpv$ then there exists a term $t(\bar{x},\bar{y})$
such that $PV$ proves
$$
\all{x}\;\all{X}\;F(\bar{x},\bar{X})\leq t(\bar{x},|\bar{X}|).
$$
\end{lemma}

Let $\lpv$ be the language which consists of function symbols for all
$FP$ functions. We denote the class of $\sbi{0}$ formulas in the language
$\lpv$ by $\sbi{0}(PV)$.

\begin{definition}
  For $\varphi(\bar{x},\bar{X})\in\sBi{0}(PV)$, $\bar{x}\in M_0^*$ and
  $\bar{X}\in M^{\boolec}$, we define $\bv{\varphi(\bar{x},\bar{X})}$
  inductively as follows:　
  \begin{itemize}
  \item if $\varphi(\bar{x})$ is an atomic formula which contains only
  number terms then
  $$
  \bv{\varphi(\bar{a})}=\left\{
                          \begin{array}{ll}
                            1&\mbox{ if }(M^*_0,M^*)\models\varphi(\bar{a}),\\
                            0&\mbox{ otherwise}.
                          \end{array}
                        \right.
  $$
  \item if $X:a\rightarrow\boolec$ then
  $$
  \bv{i\in X}=
  \left\{
  \begin{array}{ll}
    X(i)&\mbox{ if }i<a,\\
    0&\mbox{ otherwise.}
  \end{array}
  \right.
  $$
  \end{itemize}
  Let $F(\bar{x},\bar{X})\in\lpv$, $\bar{a}\in M_0^*$ and
  $\bar{A},Z\in M^{\boolec}$. We define $\bv{F(\bar{a},\bar{A})=Z}$
  by induction along the construction of $F$.
  \begin{itemize}
    \item if $F(\bar{x},\bar{X})$ is an $AC^0$ function and
    $t(\bar{x},|\bar{X}|)$ be a term bounding the length of $F(\bar{x},\bar{X})$
    as in Lemma \ref{lem:bound}. Then we have $C_0,\ldots,C_{t(\bar{x},\bar{|X|})-1}$
    such that each $C_i(\bar{x},\bar{X})$ outputs the $i$-th bit of $F(\bar{x},\bar{X})$.
    Let $\bar{a}\in M_0^*$ and $\bar{A}\in M^{\boolec}$. We define
    $$
    \bv{F(\bar{a},\bar{A})=Z}=
    \left\{
      \begin{array}{ll}
        \displaystyle\bigwedge_{0\leq i<t}\left(C_i(\bar{a},\bar{A})\leftrightarrow Z(i)\right)
        &\mbox{ if }|Z|=t,\\
        0&\mbox{ otherwise.}
      \end{array}
    \right.
    $$
  \end{itemize}
\end{definition}

What axioms the generic model satisfies depends on the complexity of the Boolean algebra.
Cook and Nguyen \cite{cn} gave a general method for constructing a minimal theory corresponding
to a given subclasses of $P$. That is, for a variety of complexity classes $C\subseteq P$ we can
give a single axiom $Ax_C$ which represents a concept for complete problems for $C$ such that
$VC=V^0+Ax_C$ captures $C$.

We expect that if the Boolean algebra represents such a computational concept then
the generic model based on that algebra satisfies $Ax_C$ and for some complexity classes we actually
have such correspondences. In this paper, we mainly treat Boolean algebras which consists of
Boolean circuits. It is possible to construct Boolean algebras for subclasses of $P$ to form
generic models for theories for such classes. We will pursue this problem in the forthcoming paper
\cite{kuroda}.

\begin{theorem}
Let $G\subseteq\mathbb{B}$ be $\mathcal{M}$-generic then $(M_0^*,M[G])\models VP$.
\end{theorem}

\noindent
(Proof). Recall that $VP=V^0+MCV$ where
$$
MCV\equiv\all a\;\all C,E\;\ext Y\;\delta_{MCV}(a,C,E,Y)
$$
and
$$
\begin{array}{ll}
\delta_{MCV}(a,C,E,Y)\equiv&\neg Y(0)\land Y(1)\land\all x<a\; x\geq 2\rightarrow\\
& Y(x)\leftrightarrow[(C(x)\land\all y<x\;(E(y,x)\rightarrow Y(y)))\\
&\phantom{Y(x)\leftrightarrow[}\lor(\neg C(x)\land\ext y<x\;(E(y,x)\land Y(y)))].
\end{array}
$$

Let $a\in M_0^*$, $C:a\rightarrow\boolec$ and $E:\pairs{a,a}\rightarrow\boolec$.
We define $Y:a\rightarrow\boolec$ by $Y(0)=0$, $Y(1)=1$ and for $x\geq 2$,
$$
Y(x)=(C(x)\land\bigwedge_{y<x}(E(y,x)\rightarrow Y(y)))\lor
(\neg C(x)\land\bigvee_{y<x}(E(y,x)\land Y(y))).
$$
Then it is readily seen that for any $C:a\rightarrow\boolec$ and $E:\pairs{a,a}\rightarrow\boolec$,
$$
(M_0,M)\models\all A\in 2^n\;eval(A,\bv{\delta_{MCV}(a,C,E,Y)})=1
$$
which implies that $\bv{\delta_{MCV}(a,C,E,Y)}\in G$.

\section{Separation problems and generic models}
\label{sec:separ-probl-gener}

Takeuti and Yasumoto relates relates the separation of complexity classes to properties of generic models.
However, their proof seems to have flaws and so we represent it with a correct proof.

\begin{theorem}[Takeuti-Yasumoto]
Let $I\subseteq\boolec$ be a $M_0$-complete ideal and $G$ be a $\mathcal{M}$-generic maximal filter over $I$.
If $P=NP$ then $(M_0^*,M[G])\models\sBi{1}$-COMP.
\end{theorem}

\noindent
(Proof). Assume that $P=NP$. Then $(M_0^*,M^*)\models P=NP$ too.
Let $\varphi(x,X,Z)\in\sBi{0}$, $t(a)$ be a term and
$\psi(x,X)\equiv\ext Z<t(|X|)\varphi(x,X,Z)$.
We will show that for any $a,b\in M_0^*$ and $X:b\rightarrow\boolec$,
$$
(M_0^*,M^*)\models
\ext Y<a\;\all x<a\;(Y(x)\leftrightarrow\ext Z<t(b)\;\psi(x,i_G(X))).
$$

Since $P=NP$, we can construct a PV function $F(x,X)$ such that
$$
\ext Z<t(|X|)\;\varphi(x,X,Z)
\rightarrow\varphi(x,X,F(x,X))\land |F(x,X)|<t(|X|)
$$
using binary search. Moreover, we can define $F(x,X)$ so that
$$
\neg\ext Z<t(|X|)\;\varphi(x,X,Z)\leftrightarrow |F(x,X)|\geq t(|X|).
$$

By the translation lemma for $PV$ formulas, we have a sequence of circuits
$$
C_{0,x}(X), \ldots, C_{t(|X|)-1,x}(X)\in\boolec
$$
such that
$$
\all X\; |X|=a\rightarrow
\all i<t(a)\;(i\in F(x,X)\leftrightarrow C_{i,x}(X)=1).
$$

Define $Z:t(b)\rightarrow\boolec$ and $Y:a\rightarrow\boolec$ by
$$
Z(i)=C_{i,x}(X)
$$
and
$$
Y(x)=\left\{
\begin{array}{ll}
\bv{\varphi(x,X,Z)}&\mbox{ if }|F(x,X)|<t(|X|),\\
0&\mbox{ otherwise.}
\end{array}
\right.
$$
Then we have
$$
\begin{array}{rcl}
Y\in G
&\Rightarrow&|F(x,X)|<t(b)\\
&\Rightarrow&\bv{\varphi(x,X,Z)}\in G\\
&\Rightarrow&
\mbox{ there is }Z:t(b)\rightarrow\boolec\mbox{ such that }
\bv{\varphi(x,X,Z)}\in G\\
&\Rightarrow&
(M_0^*,M[G])\models\ext Z<t(b)\;\varphi(x,i_G(X),Z)
\end{array}
$$
and
$$
\begin{array}{rcl}
Y\not\in G
&\Rightarrow&|F(x,X)|\geq t(b)\\
&\Rightarrow&\bv{\varphi(x,X,Z)}\not\in G\\
&\Rightarrow&
\mbox{ there is no }Z:t(b)\rightarrow\boolec\mbox{ such that }
\bv{\varphi(x,X,Z)}\in G\\
&\Rightarrow&
(M_0^*,M[G])\models\neg\ext Z<t(b)\;\varphi(x,i_G(X),Z)
\end{array}
$$

\section{Kraj\'\i\v cek forcing and T-Y forcing}
\label{sec:krajiv-cek-forcing}

In this section, we show that Kraj\'\i\v cek's forcing construction can be done by means of Takeuti-Yasumoto forcing.

\begin{theorem}
Let $(M_0,M)\models VP$ be countable and nonstandard.
If $(M_0,M)\models NP\not\subseteq P/poly$ then there exists an $M_0$-complete ideal $I\subseteq\boolec$
such that $(M_0,M[G])\models NP\not\subseteq co\mbox{-}NP$.
\end{theorem}

\noindent
(Proof). Let $(M_0,M)\models VP$ be as above. Then there exists $n\in M_0\setminus\omega$ such that
$Sat_n(X)$ is not recognized by a circuit in $(M_0,M)$, that is
$$
(M_0,M)\models\all C\in Circuit_n\;\ext X\;(Sat_n(X)\not\leftrightarrow X\models C).
$$
Let $\bar{p}=p_0,\ldots,p_{n-1}$ be a list of propositional variables in $(M_0,M)$.
Let $\varphi(X,Z)\in\sBi{0}$ be such that $\psi(X)\equiv\ext Z<|X|\varphi(X,Z)$
denotes a NP complete predicate. Define
$$
T=\{\bv{A\models X} : A:n\rightarrow\boolec\}\cup
\{\bv{\varphi(X,Z)} : Z:n\rightarrow\boolec\}.
$$
Then $T$ is consistent in $(M_0,M)$, that is, for any $T'\subseteq T$ with
$|T'|\in M_0$, it is the case that
$$
(M_0,M)\models\neg\ext P\;Prf_{EF}(\bigwedge T'\rightarrow,P).
$$
Then the theorem follows from the following observation:
\begin{lemma}
  If $T\subseteq\boolec$ is a consistent set in $(M_0,M)$ then there exists
  an $M_9$ complete ideal $I$ such that $T\subseteq G$ whenever $G$ is an
  $\mathcal{M}$-generic ultrafilter over $I$.
\end{lemma}

\noindent
(Proof). Let $T$ be consistent in $(M_0,M$ $T^*$ be the closure of $T$
under $\lor$. Note that $T^*$ is also consistent in $(M_0,M)$.
Define $I\subseteq\boolec$ by
$$
I=\{X\in\boolec : \mbox{ there exists }Z\in T^*\mbox{ such that }X<Z\}.
$$

First we claim that $I$ is an $M_0$-complete ideal. By the definition of $I$,
it is straightfoward to see that $0\in I$ and $1\not\in I$. The upward closedness
is also trivial.

Let $X_0, X_1\in I$. Then there exist $Z_0,Z_1\in T^*$ such that
$X_i<Z_i$ for $i=0,1$. Since $T^*$ is closed under $\land$, we have
$$
X_0\lor X_1<Z_0\lor Z_1\in T^*
$$
as desired. Moreover, $M_0$-completeness is proven in a similar manner.

Finally we show that if $G$ is $\mathcal{M}$-generic over $I$ then
$T\subseteq G$. To see this, let $X\in T$ and define
$$
D=\{Z\in\boolec : Z\leq X\}.
$$
We will show that $D$ is dense over $I$. Let $Y\in\boolec\setminus I$.
Then there in no $Z\in T^*$ such that $Y<Z$. Define $Z:=X\land Y$.
Then $Z\in D$ and also by the definition of $I$, we have $Z\not\in I$.

Now if $G$ is $\mathcal{M}$-generic over $I$ then there is $z\in G\cap D$.
So we have $Z\leq X$ and $X\in G$ since $G$ is a filter. \hfill$\Box$

\vspace{6pt}
Kraj\'\i\v cek gave another forcing construction.

\begin{theorem}[Kraj\'\i\v cek]\label{thm:kra2}
If $(M_0^*,M^*)\models\neg\ext P\;Prf_{EF}(\tau,P)$ for some propositional formula $\tau$. Then there
exists a $\pBi{1}$-elementary cofinal extension $(M_0^*,M')$ in which $\neg\tau$ is satisfiable.
\end{theorem}

Kraj\'\i\v cek proved this using a forcing construction with the forcing notion.
We will show that the genericl model of Theorem \ref{thm:kra2} can be obtained by
Takeuti-Yasumoto forcing.

\begin{theorem}\label{thm:kra3}
If $(M_0,M)\models\vi{1}+\neg\ext P\;Prf_{EF}(\tau,P)$ for some propositional formula $\tau$. Then there
exists a $\mathcal{M}$-generic $G\subseteq\boolec$ over some $M_0$-complete ideal such that
$(M_0,M[G])\models V^1$ and $\neg\tau$ is satisfiable.
\end{theorem}

\noindent
(Proof). We show that for any $\mathcal{G}$ which is generic in $\mathcal{P}$
there exists an $M_0$-complete ideal $I$ and an $\mathcal{M}$-generic $G$ over $I$
such that $G=\cup\mathcal{G}$. Then Kraj\'\i\v cek's theorem implies the claim.

Let $(M_0,M)$ be as above. By compactness, we construct an elementary extension $(M_0^*,M^*)$
such that there exists $c\in M_0^*$ so that
$$
(M_0^*,M^*)\models a<c
$$
for all $a\in M_0$. Construct the Boolean algebra $\boolec\subseteq M$.
Let $S\subseteq\boolec$. An $EF(S)$-proof is a sequence $P=\pairs{P_0,\ldots,P_l}$
such that each $P_i$ is either an axiom of EF, an member of $S$ or obtained from
$P_{j_0},\ldots,P_{j_k}$ for some $j_0,\ldots,j_k<i$. Note that we have
a $\sBi{0}\cup\{S\}$-formula
$$
Prf_{EF}(S,P,C)\Leftrightarrow P\mbox{ is an }EF(S)\mbox{-proof of }C.
$$
For $S\subseteq\boolec$, $C\in\boolec$ and $l\in M^*$ we say that
$S$ $l$-entails $C$ if
$$
(M_0^*,M^*)\models\ext P\; Prf_{EF}(S,P,C)\land |P|\leq l
$$
$S$ is $l$-consistent if it does not entail $0$.
Otherwise $S$ is $l$-inconsistent. Define
$$
\mathcal{P}=\{S\subseteq\boolec : S\mbox{ is }\dBi{1}
\mbox{ definable and $l$-consistent for some }l\in M_0^*\setminus M_0\}.
$$
$\mathcal{P}$ is partially ordered by the reverse inclusion.

A set $\mathcal{D}\subseteq\mathcal{P}$ is dense in $\mathcal{P}$ if for any
$S\in\mathcal{P}$ there exists $S'\in\mathcal{D}$ such that $S\subseteq S'$.
$\mathcal{D}$ is definable if
$$
\mathcal{D}=\{S\subseteq\boolec : (M_0^*,M^*)\models\eta(S)\}
$$
for some formula $\eta(S)$ with a placeholder for an unary predicate $S$.
A set $\mathcal{G}\subseteq\mathcal{P}$ is a generic if it is downward closed and
for any dense and definable set $\mathcal{D}\subseteq\mathcal{P}$,
$\mathcal{D}\cap\mathcal{G}\neq\emptyset$.

First we show that if $\mathcal{G}\subseteq\mathcal{P}$ is a generic and $G=\cup\mathcal{G}$
then $(M_0,M[G])\models\vi{1}$. Notice that the following claim can be proved in almost
the same manner as for Claim 7 of \cite{jan}.

\noindent
{\bf Claim 1}. Let $S\in\mathcal{P}$ be $l$-consistent for some $l\in M_0^*\setminus M_0$,
$\varphi(x,\bar{y},\bar{X},Z)\in\sBi{0}$, $t$ be an $\ltwo$ term and $a\in M_0$.
For any $\bar{b}\in M_0$ and $\bar{A}\in M^{\boolec}$ at least one of the following sets
is $l$-consistent:
\begin{description}
\item[(a)] $S\cup\{\neg\bv{\varphi(0,\bar{b},\bar{A},Z)} : Z:t\rightarrow\boolec\}$,
\item[(b)] $S\cup\{\bv{\varphi(a,\bar{b},\bar{A},Z_a)}\}$ for some  $Z_a:t\rightarrow\boolec$,
\item[(c)] $S\cup\{\bv{\varphi(a,\bar{b},\bar{A},Z_x)}\}\cup\{\neg\bv{\varphi(x+1,\bar{b},\bar{A},Z)} : Z:t\rightarrow\boolec\}$
for some $x<a$ and $Z_x : t\rightarrow\boolec$.
\end{description}
Note that by Claim 1, at least one of the following conditions holds for $G=\cup\mathcal{G}$;
\begin{description}
\item[(a)] $\{\neg\bv{\varphi(0,\bar{b},\bar{A},Z)} : Z:t\rightarrow\boolec\}\subseteq G$,
\item[(b)] $\{\bv{\varphi(a,\bar{b},\bar{A},Z_a)}\}\subseteq G$ for some  $Z_a:t\rightarrow\boolec$,
\item[(c)] $\{\bv{\varphi(a,\bar{b},\bar{A},Z_x)}\}\cup\{\neg\bv{\varphi(x+1,\bar{b},\bar{A},Z)} : Z:t\rightarrow\boolec\}\subseteq G$
for some $x<a$ and $Z_x : t\rightarrow\boolec$.
\end{description}
So by Forcing Theorem we have
$$
(M_0,M[G])\models IND_a(\psi(x,\bar{b},\bar{i_G(A)}))
$$
where $\psi(x,\bar{y},\bar{X})\equiv\ext z<t\varphi(x,\bar{y},\bar{X},Z)$.

Now it remains to show that

\noindent
{\bf Claim 2}. There exists an $M_0$-complete ideal such that if $\mathcal{G}$ is
$\mathcal{P}$-generic then $G=\cup\mathcal{P}$ is $\mathcal{M}$-generic over
$I$.

\noindent
(Proof of Claim).
Define
$$
\begin{array}{ll}
I=\{C\in\boolec : &\mbox{ there exists }S\in\mathcal{P}\mbox{ such that }\\
&\{X\}\cup S\mbox{ is $l$-inconsistent for some }l\in M_0\}.
\end{array}
$$
Then $I$ is an $M_0$-complete ideal.

We will show that if $\mathcal{G}$ is $\mathcal{P}$-generic then $G=\cup\mathcal{G}$
is $\mathcal{M}$-generic over $I$. To show this let $D\subset\boolec$ be a
definable and dense over $I$. It suffices to show that $G\cap(D\setminus I)\neq\emptyset$.
Define
$$
\mathcal{D}=\{S\in\mathcal{P} : S\cap(D\setminus I)\neq\emptyset\}.
$$
We claim that $\mathcal{D}$ is dense in $\mathcal{P}$.
Let $S\in\mathcal{P}$. If $S\not\in\mathcal{D}$ then $S\cap\mathcal{D}=\emptyset$.
So there exists $X\in S\setminus I$ such that $X\in D$.
We claim that $S'=\{X\}\cup S\in\mathcal{D}$. Since $S'\cap D\neq\emptyset$ is trivial,
it suffices to show that $S'\in\mathcal{P}$.

Since $X\not\in I$, $\{X\}\cup S$ is $l$-consistent for all $S$ and $l\in M_0$.
So by overspill, there exists $l\in M_0^*\setminus M_0$ such that $\{X\}\cup S$
is $l$-consistent. Thus the proof terminates. \hfill$\Box$

\section{Generic models and the pigeonhole principle}

In this section we consider the problem of whether the pigeonhole principle
holds in T-Y generic extensions. Let $F$ be an unary function symbol.
The dual weak surjective pigeonhole principle for $F$ is the following axiom:
$$
dWPHP(F)^m_n\equiv\ext Y:|Y|=m\;\all X:|X|=n\;F(X)\neq Y.
$$
where $n<m$. Jerabek \cite{jerabek} considered theories $VP$ and $V^1$ extended by axioms
$dWPHP(F)^m_n$ for all PV functions $F$. In particular,, he gave a propositional
proof system which corresponds to these theories.
\begin{definition}
The propositional proof system $WF$ extends $CF$ in the following manner;
an $WF$-proof is a sequence of circuits $X_1,\ldots,X_k$ such that each $X_i$
for $1\leq i\leq k$ is either an axiom of $CF$, obtained from $X_{j_1}\ldots,X_{j_l}$
for $j_1,\ldots,j_l<i$ or a circuit of the form
$$
\bigvee_{i<m}(r_i\not\leftrightarrow C_i(D_{i,0},\ldots,D_{i,n-1}))
$$
where $n<m$, $C_0,\ldots,C_{m-1},D_0,\ldots,D_{n-1}$ are circuits and
$r_i$ are variables which may not occur in $X_j$s for $j<i$ or $C_i$
but may occur in $D_{i,j}$s.
\end{definition}

Jerabek showed the following relations;
\begin{theorem}[Jerabek]
$V^1+dWPHP(PV)$ proves the reflection principle for $WF$.
\end{theorem}

\begin{theorem}[Jerabek]
If $V^1+dWPHP(PV)$ proves $\all\bar{x}\;\all\bar{X}\;\varphi(\bar{X},\bar{X})$
for $\varphi(\bar{X},\bar{X})$ then the propisitional translation of
$\varphi(\bar{X},\bar{X})$ have polynomial size $WF$-proofs.
\end{theorem}

Jerabek also considered a slightly weaker theory $APC_1$ which extends
$PV$ by $dWPHP(PV)$ and and established connections
between several probabilistic complexity classes. He also showed that $APC_1$ is
strong enough to manage basic combinatorial arguments such as inclusion-exclusion
principles and Chernoff's bounds.

We consider the problem of whether we can construct a generic extension of
models of $V^1+dWPHP(PV)$ in which $dWPHP(PV)$ fails. It turns out that
the answer is affirmative if we assume a condition about the complexity of
propotional proofs.
\begin{theorem}
Let $(M_0,M)\models V^1+dWPHP(PV)$ and suppose that
$$
(M_0,M)\models\all Y:|Y|=n\;
\bigwedge_{i<m}(Y(i)\leftrightarrow C_i(D_0,\ldots,D_{n-1}))\mbox{ is satisfiable.}
$$
for some $m,n\in M_0$ with $n<m$ and circuits
$C_0,\ldots,C_{m-1},D_0,\ldots,D_{n-1}\in M$. Then there exists
an $\mathcal{M}$-generic $G\subseteq\mathbb{B}_1\subseteq M$ such that
$$
(M_0,M[G])\models V^1+\neg dWPHP(F)^m_n
$$
for some PV funtion $F$.
\end{theorem}

\noindent
(Proof).
Let $(M_0,M)$ be as above and $F$ be a $PV$ function such that
$$
(M_0,M)\models\all X:|X|=n\;(|F(X)|=m\land
\all i<m\;(F(X)(i)\leftrightarrow C_i(X))).
$$
We will show that there exists an $\mathcal{M}$-generic $G\subset\mathbb{B}_1$
such that
$$
(M_0,M[G])\models\all Y:|Y|=m\;\ext X:|X|=n\;(F(X)=Y).
$$
To this end, it suffices to show that for any $Y:m\rightarrow\mathbb{B}_1$
there exists $X:n\rightarrow\mathbb{B}_1$ such that
$$
\bigwedge_{i<m}(Y(i)\leftrightarrow C_i(X))\in G.
$$

First remark that there exist $Z_0,\ldots,Z_{m-1}\in\mathbb{B}_1$ such that for all $Y\in 2^m$,
$$
\bigwedge_{i<m}(Y(i)\leftrightarrow C_i(D_0,\ldots,D_{n-1}))
$$
evaluates to true on $\pairs{Z_0(Y),\ldots,Z_{m-1}(Y)}$ in $(M_0,M)$.

For each $Y:m\rightarrow\mathbb{B}_1$, we define $X_Y:n\rightarrow\mathbb{B}_1$ as
$$
X_Y(i)=D_i(Z_0(Y),\ldots,Z_{m-1}(Y))
$$
and set
$$
S=\left\{\bigwedge_{i<m}(Y(i)\leftrightarrow C_i(X_Y)) : Y:m\rightarrow\mathbb{B}_1 \right\}.
$$
We claim that $S$ is consistent in $(M_0,M)$.

The idea for $S$ is that it forms an embedding $2^m\rightarrow 2^n$.
Moreover, for the fact that any subset $S'$ of $S$ with $|S'|\in M_0$
is consistent is guaranteed by the assumption.

Specifically, let $S'\subseteq S$ be such that $|S'|\in M_0$. By assumption we have
$$
(M_0,M)\models\all Y:|Y|=m\;\bigwedge_{i<m}(Y(i)\leftrightarrow C_i(X_Y))\mbox{ is true.}
$$
So for each $Y:m\rightarrow\mathbb{B}_1$, we have
$$
(M_0,M)\models\bigwedge_{i<m}(Y(i)\leftrightarrow C_i(X_Y))\mbox{ is a tautology.}
$$
Thus we also have
$$
(M_0,M)\models\bigwedge S'\mbox{ is a tautology.}
$$
and so
$$
(M_0,M)\models\bigwedge S'\rightarrow\mbox{ is not satisfiable.}
$$
Thus by the soundness of $EF$ in $(M_0,M)$, we have the claim.

Now, recall that any consistent set can be extended to some $\mathcal{M}$-generic
$G$ so that we have
$$
(M_0,M[G])\models\neg dWPHP(F)^m_n.
$$
\hfill$\Box$

Next we consider the problem of constructing generic extensions for $\apc$.
In fact, we can expand the Boolean algebra so that it contains
enough information order that generic extensions satisfy weak pigeonhole principles
for $PV$ functions. The idea is to construct a Boolean algebra which consists of
"randomized" circuits.

A randomized circuit with input variables $\bar{p}=p_0,\ldots,p_{n-1}$ and
random variables $\bar{z}=z_0,\ldots,z_{m-1}$ is a Boolean circuit with
variables $bar{p},\bar{z}$. We denote randomized a circuit with random varibles
$\bar{z}$ as $R\bar{z}C(\bar{p},\bar{z})$ where $C(\bar{p},\bar{z})$ is
a Boolean circuit.

\begin{definition}
$\rcircuit$ consists of circuits defined as follows:
\begin{enumerate}
\item a circuit $C(\bar{p})$ with variables $\bar{p}$ is in $\rcircuit$.
\item a randomized circuit $R\bar{z}C(\bar{p},\bar{z})$ is in $\rcircuit$.
\item if $C(x_0,\ldots,x_l)$ is a circuit wih inputs $x_0,\ldots,x_l$ and
$$
R\bar{z}C_0(\bar{p},\bar{z}),\ldots,R\bar{z}C_l(\bar{p},\bar{z})
$$
are randomized circuits in $\rcircuit$ then the circuit
$$
C(R\bar{z}C_0(\bar{p},\bar{z}),\ldots,R\bar{z}C_l(\bar{p},\bar{z}))
$$
which is obtained by replacing each $x_i$ in $C$ by $R\bar{z}C_i(\bar{p},\bar{z})$
is in $\rcircuit$.
\end{enumerate}
\end{definition}

Circuits in $\rcircuit$ may contain arbitrary number of random inputs.
If $C\in\rcircuit$ contains more than one randomized subcircuits then
we may assume that they contain the same set of random variables.

\begin{definition}
For $C\in\rcircuit$ and $A\in 2^n$ we define $eval_R(A,C)$ as follows:
\begin{enumerate}
\item if $C$ is a Boolean circuit then
$$
eval_R(A,C)=eval(A,C).
$$
\item Let $C=R\bar{z}C_0(\bar{p},\bar{z})\in\rcircuit$ We define
$$
\begin{array}{l}
|\{Z\in 2^m : eval(\pairs{A,Z},C_0(\bar{p},\bar{z}))=1\}|\succeq_03/4\cdot2^m
\Rightarrow eval_R(A,C)=1\\
|\{Z\in 2^m : eval(\pairs{A,Z},C_0(\bar{p},\bar{z}))=1\}|\preceq_01/4\cdot2^m
\Rightarrow eval_R(A,C)=0\\
\end{array}
$$
\item Let $C=C'(R\bar{z}C_0(\bar{p},\bar{z}),\ldots,R\bar{z}C_l(\bar{p},\bar{z}))$.
We define
$$
eval_R(A,C)=eval(D,C(e_0,\ldots,e_l))
$$
where $D(i)=eval_R(A,R\bar{z}C_i(\bar{p},\bar{z}))$ for $i\leq l$.
\end{enumerate}
\end{definition}

Note that we use approximate counting in the definition of the function $eval_r$
as the exact counting is not definable in the ground model.
Also note that $eval_R(A,C)$ is $\sBi{1}$-definable in $APC_1$.

\begin{definition}
For $C,C'\in\rcircuit$, we define
$$
C\leq_R C'\Leftrightarrow\all A\in 2^n(eval_R(A,C)\leq eval_R(A,C')).
$$
and
$$
C=_RC'\Leftrightarrow C\leq_R C'\land C'\leq_R C.
$$
Define $\booler=\rcircuit/=_R$.
\end{definition}

It is easy to see that $\booler$ forms a Boolean algebra. So we may define
generic extensions in the same manner as for $\boola$.
Moreover, forcing theorem for $\sBi{0}$ formulas holds for $\mathcal{M}$-generic
$G\subseteq\booler$.

We will show that $\booler$ is a suitable Boolean algebra for $APC_1$.

\begin{theorem}\label{thm:dwphp2}
If $(M_0,M)\models APC$ and $G\subseteq\booler$ is an $\mathcal{M}$-generic
over some $M_0$-complete ideal then
$$
(M_0,M[G])\models APC_1.
$$
\end{theorem}

\noindent
(Proof). Since $\boola\subseteq\booler$ it is straightforward to see that
$(M_0,M[G])\models PV$. So it suffices to show that $(M_0,M[G])\models dWPHP(PV)$.
Let $F(X)$ be a $PV$-function and $a\in M_0$ be such that $a\geq 2$.
Without loss of generality we may assume that
$$
\all X\in 2^a |F(X)|=2^{2a}.
$$
So there are circuits $C_0,\ldots,C_{2a-1}\in M$ such that
$$
\all X\in 2^a\all i<2a(F(X)(i)\leftrightarrow C_i(X)=1).
$$
Our goal is to show that
$$
(M_0,M[G])\models\ext Y\in 2^{2a}\all X\in 2^a(F(X)\neq Y)
$$
which is equivalent to the condition that there exists $Y:2a\rightarrow\booler$
such that for all $X:a\rightarrow\booler$ it is the case that
$$
\bigvee_{i<2a}(Y(i)\leftrightarrow C_i(X))\in G.
$$
We will prove the following claim which immediately implies the above assertion:

\noindent
{\bf Claim.} There exists $Y:2a\rightarrow\booler$
such that for all $X:a\rightarrow\booler$
$$
(M_0,M)\models\all A\in 2^a eval(A,\bigvee_{i<2a}(Y(i)\leftrightarrow C_i(X))).
$$

\noindent
(Proof of Claim). Let $\bar{z}=z_0,\ldots z_m$ be the list of random variables
appearing in $X:a\rightarrow\booler$. Set $Y:2a\rightarrow\booler$ to be
$Y(i)=y_i$ where $y_0,\ldots,y_{a-1}$ are flesh random variables which are
distinct from $\bar{z}$. Set
$$
P(A)=\{\pairs{\bar{b},\bar{c}} : \bar{b}\in 2^{2a},\bar{c}\in 2^a,
eval_(\pairs{A,\bar{b},\bar{c}},\bigvee_{i<2a}(y_i\not\leftrightarrow C_i(X)))=1\}.
$$
where $\bar{b}$ and $\bar{c}$ are assignments for $\bar{y}$ and $\bar{z}$ respectively.
Then it suffices to show that $(M_0,M)\models P(A)\succeq_0 3/4\cdot 2^{2a+m}$.

Let $\bar{x}=x_0,\ldots,x_{a-1}$ be flesh variables and
$$
P'=\{\pairs{\bar{b},\bar{c}} : \bar{b}\in 2^{2a},\bar{c}\in 2^a,
eval_(\pairs{\bar{b},\bar{c}},\bigvee_{i<2a}(y_i\not\leftrightarrow C_i(\bar{x})))=1\}.
$$
Then we have
$$
P'\simeq_0 2^{3a}-2^a\times|range(\bar{C})|\succeq_0 2^{3a}-2^a\cdot 2^a=2^{3a}-2^{2a}
$$
in $(M_0,M)$ where
$$
range(\bar{C})=\{Y\in 2^{2a} : \ext X\in 2^a(F(X)=Y)\}.
$$
Therefore
$$
P(A)\succeq_0(2^{3a}-2^{2a})\times 2^m=2^{2a+m}(2^a-1)\geq 3/4\cdot 2^{2a+m}
$$
for $a\geq 2$ which proves the claim. \hfill$\Box$

\begin{theorem}
Let $(M_0,M)\models\vi{1}+dWPHP(PV)$. Then there exists an $\mathcal{M}$-generic
$G\subseteq\booler$ such that $(M_0,M[G])\models\vi{1}+dWPHP(PV)$.
\end{theorem}

\noindent
(Proof). It is easy to see that the proof of Theorem \ref{thm:kra3} can be applied to
show that there exists an $\mathcal{M}$-generic $G\subseteq\booler$ such that
$(M_0,M[G])\models\vi{1}$. Then by Theorem \ref{thm:dwphp2}, we have
$(M_0,M[G])\models dWPHP(PV)$. \hfill$\Box$

\end{document}